\documentclass{article}
\usepackage{amssymb}
\usepackage{latexsym}
\usepackage[dvips]{graphicx}
\usepackage{color}

\newcommand{\boldgreek}[1]{\mbox{\boldmath$#1$}}

\newcommand{\R}{I\kern-0.37emR}

\newcommand{\ny}{n\rightarrow\infty}
\newcommand{\Q}{I\kern-0.37emP}
\newcommand{\E}{I\kern-0.37emE}
\newtheorem{THE}{Theorem}[section]

\usepackage{amsfonts}
\usepackage{amsmath}
\usepackage{color}

\newenvironment{proof}[1][Proof]{\noindent \textbf{#1.} }{\  \rule{0.5em}{0.5em}}
\begin{document}
\date{}

\title{Nonparametric Tests in Linear Model with Autoregressive Errors}
\author{\textsc{Olcay Arslan},\textsc{Ye\c{s}im G\"{u}ney}\\[1mm] \textsc{Jana Jure\v{c}kov\' a}\thanks{The research of J. Jure\v{c}kov\' a was supported by the Grant GA\v{C}R 18-01137S}, \textsc{Yetkin Tua\c{c}}\\[3mm]
\textit{University of Ankara and Charles University, Prague}}
\maketitle

\begin{abstract}
In the linear regression model with possibly autoregressive errors, we propose a family of nonparametric
tests for regression under a nuisance autoregression. The tests avoid the estimation of nuisance parameters, 
in contrast to the tests proposed in the literature. A simulation study, as well as an application of tests to real data, illustrate their good performance.\\[2mm]
\end{abstract}

\section{Introduction}
\setcounter{equation}{0}
The standard assumption of the independent and identically distributed errors in the linear regression model is often violated. Some authors admit the autoregressive structure of model errors. McKnight et
al. (2000) applied a double bootstrap method to analyze linear models
with autoregressive errors. The authors mostly estimated the regression parameters under autoregressive errors with a known innovation distribution. Alpuim and El-Shaarawi (2008) used  the ordinary least squares (OLS) estimator  under the $p$-order autoregressive (AR(p)) error
term and the normal innovations. Tua\c{c} et al. (2018) considered linear regression
model with AR(p) errors with Student's t-distribution 
and used conditional maximum likelihood estimation of 
model parameters.  In (2020), Tua\c{c} et al. proposed an
autoregressive regression procedure based on the skew-normal and skew-$t$
distributions. G\"{u}ney et al. (2020a) considered the conditional
maximum Lq-likelihood (CMLq) estimation method for the autoregressive error
terms regression models under normality assumption.

In the real life applications,  
the data sets may contain outliers and their distribution can be heavy-tailed. Then we should take recourse to nonparametric models without a specific distribution assumptions. The most powerful tools for estimation and other inference in this area are the regression and autoregression quantiles. However, the most convenient for testing are the ranks of observations or the ranks of residuals and their extensions, as the regression and autoregression rank scores. The quantile regression was 
introduced by Koenker and Bassett (1978) and by their followers. It is an approach to model the conditional quantile function of
response variable depending on covariates. The regression rank scores were introduced by Gutenbrunner and Jure\v{c}kov\' a (1992) and a generalization of this concept to the autoregressive model is due to Koul and Saleh (1995). The class of regression rank scores tests was developed by Gutenbrunner and Jure\v{c}kov\' a (1992) and Gutenbrunner et al. (1992). The optimal autoregression rank scores tests in the AR model were constructed by Hallin and Jure\v{c}kov\' a (1999). El Bantli and Hallin (2001) constructed the Kolmogorov-Smirnov type test in AR model based on the autoregression rank scores, following the KS test in the linear model test by Jure\v{c}kov\' a (1991). The averaged autoregression quantiles and their asymptotics were studied by G\"{u}ney et al. (2020b).
{In the present paper, we assume that our data follow a linear regression model whose model errors can be possibly autoregressive. The corresponding  probability distributions are generally unknown, only satisfy some assumptions. In this setup, we shall verify
 the hypothesis of no linear regression under possible nuisance autoregression of model errors.
 
 \section{Statement of the model}
\setcounter{equation}{0}
We consider the linear regression model of order $s,$ whose model errors follow a stationary autoregressive  process of order $p:$
\begin{eqnarray}
\label{1}
&&y_{t}=\beta_0+\mathbf x_t^{\top}\boldgreek\beta^{\ast}+\varepsilon _{t}=\beta_0+x_{t1}\beta _{1}+...+x_{tp}\beta _{p}+\varepsilon _{t}, \\
\label{1a}   
&&\varepsilon _{t}=\varphi_0+\varphi _{1}\varepsilon_{t-1}+u_{t}+\ldots+\varphi _{p}\varepsilon_{t-p}, \; t=1,2,...,n,\\ 
&&\boldgreek\beta^{\ast}=(\beta_1,\ldots,\beta_s)^{\top}.\nonumber 
\end{eqnarray}%
Here $y_{t}$ is the response variable, $\mathbf x_t=(x_{t0},\ldots,x_{ts})^{\top}$ are the regressors and $\beta_j, \; j=0,\ldots,s$
 are unknown regression parameters. We assume that $x_{t0}=1$ for all $t,$ hence that $\beta_0$ is an intercept.
 
Moreover, $\varphi _0,\varphi_1,\ldots,\varphi_p$ are unknown
autoregression parameters, where the intercept $\varphi_0$ is added for mathematical convenience and can be 0. The innovations $u_{t}$ are assumed being independently
and identically distributed (\textit{i.i.d.}) with a continuous distribution
function $F$ and density $f,$ generally unknown but satisfying
\begin{equation}\label{2.2}  
E\left( u_{t}\right) =0, \quad Var\left(u_{t}\right) =\sigma ^{2}<\infty.
\end{equation}
 The stationarity condition
requires that all roots of the equation $z^{2}-\varphi _{1}z^{1}=0 -\varphi_{2}z^{p-2}-...-\varphi _{p}z^{t-p}=0
$ are inside the unit circle (Brockwell and
Davis, 1987). Note that the error model given in equation (\ref{1a}) is a strictly
stationary process, and so the $\varepsilon _{t}$ share a common marginal distribution and thus share the same quantiles. The intercept term $\varphi_0$ in model (\ref{1}) is included for the identifiability of the autoregression quantiles, and can be equal to 0.} {The distribution function $F$ of $u_{t}$ is unknown}, but we assume that it
 is increasing on the set $\{u:0<F(u)<1\}$. 
{Because of  the identifiability, } 
we assume that the starting observations $(y_{-p+1},\ldots,y_0)$ are known. 

For the convenience, we also write (\ref{1a}) in the form          
\begin{equation}\label{1b}
u_{t}=\Phi \left( B\right) \varepsilon_{t}
\end{equation}
where $B$ is called the backshift operator. Then the linear
regression model with AR(p) error term given in equation (\ref{1}) can be expressed as
\begin{eqnarray}\label{4}
\Phi \left( B\right) y_{t}=y_{t}-\varphi_0-\varphi _{1}y_{t-1}+...\varphi _{p}y_{t-p},\\
\mathbf{\Phi }\left( B\right) \mathbf{x}_{t}=\mathbf{x}_{t}-\varphi _{0}-\varphi _{1}\mathbf{x}_{t-1}
+\ldots\varphi _{p}\mathbf{x}_{t-p},\nonumber  
\end{eqnarray} 
hence
\begin{equation}\label{3}
\Phi \left( B\right) y_{t}=\left( \mathbf{\Phi }\left( B\right) \mathbf{x}_{t}\right)^{T}\boldgreek\beta+u_{t}.  
\end{equation}%
In the model (\ref{1}), we construct the tests of the hypothesis:
\begin{description}
\item $\mathbf H_0: \; \boldgreek\beta^{\ast}=\mathbf 0, \; \mbox{ with } \; \beta_0,~ (\varphi_0, \varphi_1,\ldots,\varphi_p)^{\top}\neq \mathbf 0 \; \mbox{ unspecified.}$                     \end{description}
Our tests are nonparametric; the test of $\mathbf H_0$ is based on the autoregression rank scores and on the linear autoregression rank statistics for the model
(\ref{1}), (\ref{1a}) without regression.

\section{Rank tests for $\mathbf H_0$}
\setcounter{equation}{0}
We shall be testing the absence of regression
$$\mathbf H_0: \; \boldgreek\beta^{\ast}=\mathbf 0, \qquad  \beta_0,  \; \varphi_0, \; \varphi_1,\ldots,\varphi_p \; \mbox{unspecified}.$$
The usual alternative is the local (Pitman) regression
\begin{equation}\label{Pitman}
\mathbf K_{n}: \; \boldgreek\beta^{\ast}=\boldgreek\beta_n^{\ast}=n^{-1/2}\boldgreek\beta_{x}^{\ast}, \; \mbox{ with } \; \boldgreek\beta_{x}^{\ast}\in\R^p \; \mbox{ fixed.}
\end{equation}
Under $\mathbf H_0,$ the observations follow the model
$$y_t=\beta_0+\varepsilon_t, \; t=1,\ldots,n.$$
 The hypothesis $\mathbf H_0$ in fact means the hypothetical autoregressive model
 \begin{equation}\label{201}
\mathbf H_0: \; y_t= \beta_0+\varphi_0+\varphi_1 y_{t-1}+\varphi _{1}y_{t-1}+\ldots+\varphi_{p}y_{t-p}+u_{t}, \; t=1,\ldots,n
\end{equation}
which we like to test against the alternative $\mathbf K_{n}.$
 Let
\begin{equation}\label{201a}
{\mathbf y}_t^* = (y_t,y_{t-1},\ldots,y_{t-p})^{\top}, \;
{\mathbf y}_t= (1,\,y_t, y_{t-1},\ldots,y_{t-p})^{\top}, \; t=0,\ldots,n-1 
\end{equation}
and consider the random matrices of the respective orders $n\times p$ and $n\times (p+1)$ 
\begin{eqnarray}\label{201b}
 & \mathbf Y_n^{\ast}=\left[\begin{array}{c}\mathbf y_1^{\ast\top}\\ 
                                             \cdots\\
                                             \mathbf y_n^{\ast\top}\\                                                                     
                                             \end{array}\right], \qquad 
   \mathbf Y_n=\left[\begin{array}{c}\mathbf y_1^{\top}\\ 
                                             \cdots\\
                                             \mathbf y_n^{\top}\\                                                                     
                                             \end{array}\right]. &                                          
 \end{eqnarray} 
 For convenience, denote also 
\begin{eqnarray*}
 &\mathbf x_t^{\ast}=(x_{t1},\ldots,x_{tp})^{\top}, \quad \mathbf x_t=(1, x_{t1},\ldots,x_{ts})^{\top}=(1, \mathbf x_t^{\ast \top})^{\top}& \\ 
 & \mathbf X_n^{\ast}=\left[\begin{array}{c}\mathbf x_1^{\ast\top}\\ 
                                             \cdots\\
                                             \mathbf x_n^{\ast\top}\\                                                                     
                                             \end{array}\right], \qquad 
   \mathbf X_n=\left[\begin{array}{c}\mathbf x_1^{\top}\\ 
                                             \cdots\\
                                             \mathbf x_n^{\top}\\                                                                     
                                             \end{array}\right]. &                                          
 \end{eqnarray*}

The autoregression rank scores $\widehat{\mathbf a}_n(\alpha)=\left(\hat{a}_{n1}(\alpha),\ldots,\hat{a}_{nn}(\alpha)\right)^{\top}$ 
under hypothesis $\mathbf H_0$ are defined as the  solution vector of the linear programming problem
\begin{equation} \label{RR}
\left\{
\begin{array}{lll}
       \sum_{t=1}^n y_t\hat{a}_{nt}(\alpha) :& = &\mbox{ max }\\[2mm] 
	     \sum_{t=1}^n\left(\hat{a}_{nt}(\alpha)-(1-\alpha)\right)&=&0\\[2mm]   
        \mathbf Y_n^{\ast\top}\left(\widehat{\mathbf a}_n(\alpha)-(1-\alpha){\mathbf 1}_n \right)& =& \mathbf 0\\[2mm]
 \widehat{\mathbf a}_n(\alpha) \in [0,1]^n,\quad 0\leq \alpha \leq 1.&&\\
\end{array}
\right.
\end{equation}
The autoregression rank scores are  
{\em autoregression-invariant}. More precisely, (\ref{RR}) implies that 
$\widehat{{\mathbf a}}_{n}(\alpha)$   can be also formally written as a solution of the linear program
$$ \left\{
\begin{array}{lll}
      \sum_{t=1}^n u_t\hat{a}_{nt}(\alpha) :& = &\mbox{ max }\\[2mm]   
	    \sum_{t=1}^n\left(\hat{a}_{nt}(\alpha)-(1-\alpha)\right)&=&0\\[2mm] 
	    \mathbf Y_n^{\ast\top}\left(\widehat{\mathbf a}_n(\alpha)-(1-\alpha){\mathbf 1}_n \right)& =& \mathbf 0\\[2mm]
      \widehat{\mathbf a}_n(\alpha) \in [0,1]^n,\quad 0\leq \alpha \leq 1.&&\\
\end{array}
\right.   
$$	 
where $\mathbf u_n  =(u_1,\ldots ,u_n)^{\top}$ is the unobservable white noise process.

We shall construct a family of new tests of the hypothesis $\mathbf H_0$ for the model (\ref{1}), based on autoregression rank scores, and analyze the asymptotic distribution of the test criterion under the null hypothesis as well as under contiguous alternatives. Surprisingly,  
no preliminary estimation of $\boldgreek\varphi$ is needed in order to compute autoregression rank score statistics, in contrast with the aligned rank methods
(Puri and Sen   and others).

The (unknown) density $f$ of $u_t$ is assumed to belong to the family ${\mathcal F}$ of exponentially tailed densities,  satisfying (\ref{2.2})
and the following conditions on the tails:
\begin{enumerate}
\item[{\bf (F1)}] $f$ is  positive and absolutely continuous,
with  a.e. derivative $f^{\prime}$ and finite Fisher information
$\mathcal I (f) =
\int \left(\frac{f^{\prime}(x)}{f(x)}\right)^2f(x){\rm d}x<\infty$;
moreover,
there exists $K_f\geq 0$ such that $f$ has two bounded derivatives $f^{\prime}$ and $f^{\prime \prime}$ for all $\vert x \vert > K_f$;
\item[{\bf (F2)}] $f$ is monotonically decreasing to $0$ as $x\rightarrow
 \pm \infty$  and 
$$\lim_{x\rightarrow -\infty}\frac{-\log F(x)}{b\vert x\vert^{r}}=
\lim_{x\rightarrow\infty}\frac{-\log(1-F(x))}{b\vert x\vert^{r}}=1
$$
for some $b>0$ and $r\geq 1$.
\end{enumerate}
Other properties of densities in $\mathcal F$ are summarized in \cite{Hallin1999}.\\[2mm]

Moreover, we impose the following conditions on the regression matrix $\mathbf X_n^{\ast}:$
\begin{enumerate}
\item[{\bf (X1)}] The matrix $\mathbf A_n=n^{-1}\sum_{t=1}^n\mathbf X_n^{\ast\top}\mathbf X_n^{\ast}$ is positive definite of order $s$ for $n\geq n_0.$
\item[{\bf (X2)}] $n^{-1}\sum_{t=1}^{n}\|\mathbf x_{nt}\|^4=O(1)$ as $\ny.$
\item[{\bf (X3)}] $\lim_{\ny}\max_{1\leq t\leq n}\left\{n^{-1}\mathbf x_{nt}^{\top}\mathbf A_n^{-1}\mathbf x_{nt}\right\}=\mathbf 0.$ 
\end{enumerate}
Define $\mathbf D_n=n^{-1}\mathbf Y_n^{\top}\mathbf Y_n$ and
 \begin{equation}\label{2.3}
 \mathbf H_n=\mathbf Y_n^{\top}(\mathbf Y_n^{\top}\mathbf Y_n)^{-1}\mathbf Y_n, \qquad \widehat{\mathbf X}_n^{\ast}= \mathbf H_n\mathbf X_n^{\ast}
\end{equation} 
the projection matrix and the projection of $\mathbf X_n^{\ast}$ on the space spanned by the columns of $\mathbf Y_n$, respectively. Moreover, let  
\begin{equation}\label{2.3a}
  \mathbf Q_n=
 n^{-1}(\mathbf X_n^{\ast}-\widehat{\mathbf X}_n^{\ast})^{\top}(\mathbf X_n^{\ast}-\widehat{\mathbf X}_n^{\ast}).
\end{equation}
The (random) matrices $\mathbf D_n$ and $\mathbf Q_n$ are  of the respective orders $(p+1)\times (p+1)$ and $s\times s$.
  We shall assume that
\begin{equation}\label{2.3b}	
 \lim_{\ny}\mathbf D_n=\mathbf D, \qquad \lim_{\ny}\E(\mathbf Q_n)=\mathbf Q 
\end{equation}
where $\mathbf D$ and $\mathbf Q$ are positive definite matrices.
 
Choose a nondecreasing, square integrable score generating function $J : (0,1)\rightarrow \R$, such that $J(1-u) =-J(u),\; 0 < u < 1$
and that $J^{\prime}(u)$ exists for
$u\in (0,\, \alpha_0)\cup (1-\alpha_0,\, 1)$
 and, in this domain,  
satisfies the  Chernoff-Savage condition
\begin{equation}
\label{chernoffsavage}
|J^{\prime}(u)|\leq c(u(1-u))^{-1-\delta}, \ 0<\delta<\frac{1}{4}.
\end{equation}
Define the scores generated by $J$ as $\widehat{\mathbf b}_n=(\hat{b}_{n1},\ldots,\hat{b}_{nn})^{\top}$ with
 \begin{equation}\label{3.1}
 \hat{b}_{nt} 
= -\int_{0}^1 J(u) {\rm d} \hat a_{nt} (u), \quad t = 1, \ldots, n.	
 \end{equation}
The proposed tests of $\mathbf H_0$ are based on the linear autoregression rank statistics of the form
  \begin{equation}\label{3.3}
	{\mathbf S}_{n} =  n^{-\frac{1}{2}}(\mathbf X_n^{\ast}-\widehat{\mathbf X}_n^{\ast})^{\top}\widehat{\mathbf b}_{n},	
\end{equation}
and  for testing $\mathbf H_0$ against $\mathbf K_{n}$ we propose the criterion
\begin{equation}\label{2.9}
\mathcal T_n=\mathbf S_n^{\top}\mathbf Q_n^{-1}\mathbf S_n/A^2(J)
\end{equation}
where
\begin{equation}\label{2.10}
A^2(J)=\int_0^1 (J(t)-\bar{J})^2{\rm d} t, \qquad \bar{J}=\int_0^1 J(t){\rm d} t.
\end{equation}
The typical choices of $J$ are:
\begin{description}
	\item{(i)} 
Wilcoxon scores (optimal for $f$ logistic) : $J(u) = u - \frac 12, \; 0 < u < 1.$ The scores are $\hat b_{n;t}= - \int_0^1 (J(u)-\frac 12)d\hat{a}_t(u) =
\int_0^1 J^2(u)du-\frac 12$ while $A^2(J) = \frac{1}{12}$ and $\gamma(J, F) = \int f^2(x)dx.$ 
\item{(ii)} 
 Normal (van der Waerden) scores (asymptotically optimal for $f$ normal): $J(t) = \Phi^{-1}(u), \; 0 < u < 1, \Phi$ being the d.f.
of standard normal distribution. Here $A^2(J) = 1$ and $\gamma(J, F) =
\int f(F^{-1}(J(x)))dx$. 
\item{(iii)} Median (sign) scores: $J(u) = \frac 12~sign(u - \frac 12), \; 0 < u < 1.$ 
\end{description}
Notice that the test statistic $\mathcal T_n$ requires no estimation of nuisance parameters, since the functional $A(J)$ depends only on the score function and not on (the unknown) $F.$  We shall show that the asymptotic distribution of $\mathcal T_n$ under $\mathbf H_0$ is central $\chi^2$ with $s$ degrees of freedom, hence it is asymptotically distribution free. Under $\mathbf K_{n}$ it is noncentral $\chi^2$ with $s$ degrees of freedom and noncentrality parameter dependent on $J$ and $F$ but not on the nuisance parameters. In this way, it is asymptotically equivalent to the rank test of $\mathbf H_0$ in the situation without nuisance autoregression. 
\section{ Asymptotic behavior of the test of $\mathbf H_0$}
\setcounter{equation}{0}

Let us return to the model (\ref{1}). Assume that the matrices $\mathbf X_n^{\ast}$ and $\mathbf Y_n^{\ast}$ satisfy
conditions (\ref{2.3})--(\ref{2.3b}). We want to test the hypothesis
$$\mathbf H_0: \; \boldgreek\beta^{\ast} = \mathbf 0 \quad (\beta_0, \; \boldgreek\varphi \; \mbox{ unspecified })$$ against the alternative 
$$\mathbf K_n: \; \boldgreek\beta^{\ast} = 
n^{-1/2}\boldgreek\beta_{x}^{\ast} \quad  (\boldgreek\beta_{x}^{\ast}\in\R_s \; \mbox{ fixed }).$$
Let $\widehat{\mathbf a}_n(\alpha)=(\hat{a}_{n1}(\alpha),\ldots,\hat{a}_{nn}(\alpha))$ be the autoregression rank scores corresponding to the submodel under 
$\mathbf H_0,$ i.e.
$$y_t= \beta_0+\varphi_0+\varphi_1 y_{t-1}+\varphi _{1}y_{t-1}+\ldots+\varphi_{p}y_{t-p}+u_{t}, \; t=1,\ldots,n.$$
Let $J : (0,1)\mapsto \R$ be a nondecreasing and square integrable score-generating
function such that $J(1-u) =-J(u), \; 0 < u < 1$, satisfying (\ref{chernoffsavage}). Define the scores $\widehat{\mathbf b}_n=(\hat{b}_{n1},\ldots,\hat{b}_{nn})^{\top}$ by the relation (\ref{3.1}). Consider
the test statistics
$\mathcal T_n=\mathbf S_n^{\top}\mathbf Q_n^{-1}\mathbf S_n/A^2(J)$
defined in (\ref{3.3})--(\ref{2.10}).
The test is based on the asymptotic distribution of $\mathcal T_n$ under $\mathbf H_O$, described in the following theorem.
\begin{THE}\label{Theorem1}
Assume that the distribution $F$ of the innovations $u_t$ satisfies (F1)--(F4) and 
the regression matrix $\mathbf X_n$ satisfies (X1)--(X4). Let $\mathbf T_n$ be generated by the function $J$ satisfying (\ref{chernoffsavage}), nondecreasing and square integrable on $(0,1).$ 
\begin{description}
	\item{(i)} Then, under $\mathbf H_0,$ the asymptotic distribution of $\mathcal T_n$ is central $\chi^2$ with $s$ degrees of freedom.
	\item{(ii)} Under $\mathbf K_n,$ the asymptotic distribution of $\mathcal T_n$ is noncentral $\chi^2$ with $s$ degrees of freedom and the noncentrality parameter
\begin{eqnarray}\label{5.4}	
&&\eta^2=\boldgreek\beta_x^{\top}\mathbf Q\boldgreek\beta_x\cdot\gamma^2(J,F)/A^2(J)\nonumber\\
\mbox{ where }&& \\
&&\gamma(J,F)= -\int_0^1 J(v){\rm d}f(F^{-1}(v)).\nonumber
\end{eqnarray}
\end{description}
\end{THE} 
Hence, the test rejects $\mathbf H_0$ on the significance level $\tau\in(0,1)$ if
$\mathcal T_n>\chi_s^2(1-\tau)$ where $\chi_s^2(1-\tau)$ is the $100(1-\tau)\%$-quantile of the $\chi^2$ distribution with $s$ degrees of freedom. The asymptotic distribution under $\mathbf K_n$ also shows that the Pitman efficiency
of the test coincides with that of the classical rank test in the situation without the autoregressive errors.\\[3mm]
\begin{proof} 
\begin{description}
  \item{(i)} It follows from \cite{Hallin1999}, Theorem 3.3 and \cite{GJKP}, Theorem 4.1, that under $\mathbf H_0,$ the linear autoregression rank scores statistic admits the representation 
$$\mathbf Q_n^{-1/2}\mathbf S_n=n^{-1/2}\mathbf Q_n^{-1/2}(\mathbf X_n^{\ast}-\overline{\mathbf X}_n^{\ast})^{\top}\widetilde{\mathbf b}_n+o_p(1)$$
as $\ny,$ where $\widetilde{\mathbf b}_n=(\tilde{b}_{n1},\ldots,\tilde{b}_{nn})^{\top}$
and $\tilde{b}_{nt}=J(F(u_{nt})), \; t=1,\ldots,n.$ Then (i) follows from the central limit theorem.
  \item{(ii)} The same representation holds also under the sequence of alternatives $\mathbf K_n,$ which is contiguous with respect to the sequence of null distributions with the densities $\prod_{t=1}^nf(u_t).$ 
  \end{description}
\hfill\end{proof}

\end{document}